\documentclass[11pt]{article}      
\usepackage{amsmath,amsthm,amsfonts,graphicx}
\usepackage{multirow}
\usepackage{multirow}
\usepackage{xcolor}
\def\smallddots{\mathinner{\raise7pt\hbox{.}\raise4pt\hbox{.}\raise1pt\hbox{.}}} 
\def\smallsdots{\mathinner{\raise1pt\hbox{.}\raise4pt\hbox{.}\raise7pt\hbox{.}}}

\DeclareMathOperator{\diag}{diag}

\DeclareMathOperator{\rank}{rank}

\newtheorem{theorem}{Theorem}[section]

\numberwithin{equation}{section}
\numberwithin{table}{section}
\newtheorem{lemma}{Lemma}[section]

\newtheorem{corollary}{Corollary}[section]

\newtheorem{algorithm}{Algorithm}[section]
\newtheorem{example}{Example}[section]

\newtheorem{remark}{Remark}[section]

\begin{document}
 
\title{Iterative Refinement and Oversampling \\ for  Low Rank Approximation} 
\author{Victor Y. Pan} 
\author{Victor Y. Pan$^{[1, 2],[a]}$, Qi Luan$^{[1],[b]}$, and  Soo Go$^{[1],c]}$
\\ \\
$^{[1]}$ Ph.D. Programs in  Computer Science and Mathematics \\
The Graduate Center of the City University of New York \\
New York, NY 10036 USA \\
$^{[2]}$ Department of Computer Science \\
Lehman College of the City University of New York \\
Bronx, NY 10468 USA \\
$^{[a]}$ victor.pan@lehman.cuny.edu \\ 
http://comet.lehman.cuny.edu/vpan/  \\
$^{[b]}$
 qi\_luan@yahoo.com \\
$^{[c]}$ sgo@gradcenter.cuny.edu
\\
} 
\date{}

\maketitle

 
\begin{abstract}  
Iterative refinement is particularly popular for numerical solution of linear systems of equations. We  extend it to Low Rank Approximation of a matrix  (LRA) and observe close link
of the resulting algorithm to oversampling techniques, commonly used in randomized LRA algorithms.
 We elaborate upon this link and revisit oversampling and  some  efficient   randomized LRA algorithms. Applied with sparse sketch matrices
 they run significantly faster and in particular yield Very Low Rank Approximation (VLRA) at sublinear cost,
 using much fewer 
 scalars and flops than the input matrix has entries. This is achieved at
  the price of  deterioration of output accuracy, but according to our  formal  
 and empirical study subsequent  oversampling improves accuracy  to near-optimal level under the spectral norm for a large sub-class of  matrices with fast decaying spectra of singular values. 
\end{abstract} 
 




\bigskip
\medskip


\section{Introduction:
 LRA  of a matrix,  its iterative refinement, link to oversampling, and sparse sketches}\label{sintr} 


\subsection{Brief overview}\label{sbrf}
Classical  {\em iterative refinement}
(see pointers to the bibliography in Sec. \ref{srelwrk}) is particularly popular for numerical solution of linear systems of 
equations. We naturally extend it to {\em Low Rank Approximation of a matrix  (LRA)},
observe its close link to LRA with oversampling, and analyze both refinement and oversampling.
Here is a brief outline of our study.

 Let an LRA algorithm $\mathbb A=\mathbb A_r$ be applied to an $m\times n$ matrix $E_0:=M$ and let it output a rank-$r$ approximation $Y_0:=X_0:=\mathbb A(E_0)$ of $M$. \footnote{Here and hereafter 
``rank-$r$ approximation" means ``rank at most $r$ approximation".} Apply the algorithm $\mathbb A$ to the error matrix 
$E_1=E_0-X_0$ and add 
$X_1:=\mathbb A(E_1)$
to $X_0$ to obtain improved approximation
 $Y_1:=X_0+X_1$ to $M$. By repeating this step recursively $h$ times, compute the matrices $X_i:=\mathbb A(E_i)$, $Y_i:=\mathbb A(Y_{i-1}+X_i)$, and $E_{i+1}:=E_i-X_i$,  for $i=1,2,\dots,h$.
 Finally   output the matrix $\mathbb A(Y_h)$ as a   rank-$r$  approximation of $M$. 
The  {\em compression} of higher rank  approximation of a matrix  into its lower
 rank approximation  also occurs in the
 customary random {\em oversampling}  for LRA.

We revisit random sampling  LRA algorithms
with oversampling,  accelerate them by using sparse sketch matrices and
in particular yield {\em Very Low Rank Approximation (VLRA)} at sublinear cost. We achieve this
at the price of some deterioration of the  accuracy of the output LRA, but according to our  formal  
 and empirical study, subsequent  oversampling improves accuracy  to near-optimal under the spectral norm for a large sub-class of the matrices with fast decaying spectra of singular values.
 So far, our study has showed no additional benefits for LRA from  application of iterative refinement versus 
 standalone oversampling.

 
\subsection{The State of the Art of LRA}\label{sstart}

LRA is a hot topic of Numerical Linear 
Algebra. We first briefly recall the State of the Art,  referring  to  \cite{HMT11,TYUC17,MT20} for broader expositions, and then outline our results for LRA.

{\bf 1.1.1. Representation and the power of LRA.}  
Assume that a {\em black box matrix}  $M\in \mathbb C^{m\times n}$ is given by an oracle (black box subroutine) for computing its product with a vector rather than by its entries; then linear operations with  matrices  
 can be reduced to   linear operations with vectors.\footnote{The   sketch  of a matrix in \cite{W14,TYUC17}  is its product with a matrix of smaller size, that is,  with a  small set of vectors.} In this way  we can handle
 matrices of immense size that come from the study of  Big Data and are too large to be accessed otherwise (see further motivations in \cite[Sec. 1.3]{TYUC17}). Furthermore,  quite typically
 the matrices representing Big Data are close to low rank matrices  \cite{UT19}, that is,  admit LRA, with which we can
 operate efficiently, even where $M$ has enormous size. 
  
 Now assume that $m\ge n$ (otherwise we could have worked with the $n\times m$ transpose $M^T$) and  that we are given a positive target rank $r< n$ or  
a fixed target range for such a rank, although actually a proper choice of such targets is a challenge (see Remark \ref{reovrsmpl}).  Then define LRA of $M$ by a triplet  $\{X,Y,Z\}$ such that
\begin{equation}\label{eqXYZ}
X\in \mathbb C^{m\times k}, Y\in \mathbb C^{k\times \ell},~Z\in \mathbb C^{\ell\times n},
\end{equation}
 \begin{equation}\label{eqklmn}
 k\le m,~\ell\le n,~
\min\{k,\ell\}= r, 
\end{equation}
 and the spectral and/or Frobenius norm
 of the matrix
 $M-XYZ$ is small in context, say, close to optimal.



{\bf 1.1.2. LRA from top SVD; SVD   of an LRA.}
By setting to 0 all  singular values of $M$ except for its $r$ top  (largest) ones, turn SVD of $M$ into its {\em $r$-top SVD} and turn $M$ into its $r$-{\em truncation}  $M_r$, of rank $r$. 
 By virtue of the  Eckart-Young-Mirsky  theorem \cite[Thm. 2.4.8]{GL13}, $r$-top SVD defines an optimal rank-$r$ 
approximation
of $M$ under  both
 spectral and Frobenius norms, that is, $X=M_r$ minimizes these norms of $M-X$ over all matrices $X$ of rank $r$.  

Computation of SVD of  $M$    involves  memory space for
more than $mn$ scalars and order of $mn^2$ floating point operations  (flops)   (see \cite[Fig. 8.6.1]{GL13}). This is  too
expensive for 
 LRA of a large size matrix, while conversely,
given a triplet $\{X,Y,Z\}$  satisfying (\ref{eqXYZ}) and (\ref{eqklmn}), we can  compute  SVD of $XYZ$ 
 by using just $(m+k)k+
(n+\ell)\ell$  scalars and $O(mk^2+n\ell^2)$ flops
 (see Appendix  \ref{slrasvd}).

{\bf 1.1.3. LRA via random sampling and oversampling.} Random sampling with oversampling enables  fast computation of close LRA  (see \cite{CW09,M11,HMT11,W14,TYUC17,N20}, and the  references therein). 
For  any $M\in \mathbb C^{m\times n}$
and target rank $k<n\le m$, random sampling  {\em HMT} algorithm of \cite{HMT11}  
 successively generates  a random  {\em sketch matrix} $H\in\mathbb C^{n\times (k+p)}$
for some positive {\em oversampling integer $p$}, usually a small constant,   
computes the product $MH$, and  computes rank-$k$ approximation of $M$. 
\begin{remark}\label{reovrsmpl}
``In practice, the target rank $k$ is rarely known in advance.
Randomized algorithms
 are usually implemented in an adaptive fashion where the number of samples is increased until the error norm satisfies the desired tolerance"
(see \cite[Sec. 4.2]{HMT11}).
\end{remark} 

In the best studied case  the sketch matrix $H$ is random Gaussian, filled with independent standard Gaussian (normal) random variables. In that case an expected Frobenius error norm of the output LRA is within  a factor of $\epsilon=(1+\frac{k}{p-1})^{1/2}$ from the optimal value $||M-M_r||_F$  \cite[Thm. 10.5]{HMT11}, and 
 the algorithm uses orders of $(m+n)(k+p)$ scalars and of $(k+p)mn$ flops.
  The integer $p$ and the size of the   matrix  $H$   explode as $\epsilon \mapsto 0$, but  LRA for  $\epsilon$ as large as 1/2 or 1 is still valuable.

 
  In other LRA algorithms (e.g., in \cite{DMM08,CW09,TYUC17,N20})
 $M$ is both pre- and post-multiplied by  random  sketch  matrices $F\in \mathbb C^{k\times m}$ and $H\in \mathbb C^{n\times \ell}$, respectively.
In particular,
 the {\em TYUC} algorithm of \cite{TYUC17},  modifying those of \cite[Thm. 4.7]{CW09}
and \cite[Thm. 4.3, display 1]{CW09}, and the
generalized Nystr{\"o}m approximation {\em GN} of \cite{N20} first compute  the products $FM$ and $MH$ for random  sketch matrices $F$ and $H$
for $k\le \ell$, say, for $\ell =2k$,  
and then compute rank-$k$ approximation of $M$.
In the case of random  Gaussian  
 sketch matrices $F$ and $H$ the expected Frobenius error norm  of the output LRA is within a factor of 2 from optimal (namely, it does not exceed the  optimal one for rank-$\frac{k}{2}$ approximation).                                                                                                                                                                                                                                                                                                                                                                                                                                                                                                                                                                                                                                                                                                                                                                                                                                                                                                                                                                                                                                                                   
For $\ell=2k$ the TYUC and GN algorithms use about  $(m+n)\ell$ scalars
 (or slightly less in  \cite{N20}) and  $O(\ell mn)$ 
 flops. 
 
  The computation becomes  a little more  costly  for supporting the same 
  LRA error bounds under the spectral rather than  Frobenius norm  
 (see \cite[Secs. 10 and 11]{HMT11},
 \cite[Ch. 6.2]{W14}, \cite[Thm. 3.4]{G15}).

 \subsection{Acceleration of randomized LRA  by using sparse sketches}\label{sprgr}

Hereafter  $\sigma_j(M)$ denotes the $j$th largest 
singular value of an
$m\times n$ matrix $M$, $j=1,2,\dots$; $\sigma_j(M)=0$ for $j>\rank(M)$, and let 
$M$ be an $(r,\rho)$-{\em matrix}, that is, let it have gaps in both segments $[\sigma_1(M),\sigma_{r+1}(M)]$ and
$[\sigma_{r+1}(M),\sigma_{\rho+1}(M)]$ of its singular values for $r<\rho\le n\le m$.
We seek rank-$r$ approximation of
$(r,\rho)$-{\em matrices} $M$.

We rely on combining {\em three simple observations}, namely, that (i)  random sampling with oversampling is the bottleneck  stage of HMT, TYUC, and GN algorithms, (ii) which we can greatly accelerate by using sparse   sketch matrices (in particular abridged SRHT matrices of Appendix \ref{spreprmlt}),   and
(iii) while   
the accuracy of the LRA   output in step (ii) tends to deteriorate with sparse sketch matrices, we  readily improve it  to a near-optimal level for a large subclass of $(r,\rho)$-matrices $M$. 
\medskip

Our LRA
Alg. \ref{algesc1} of Sec. \ref{sspdsc} relies on these observations. 
It first computes (1) a crude rank-$\rho$ approximation $M(\rho)$  of $M$;
  then (2) $r$-top SVD of $M(\rho)$, output as rank-$r$ approximation of $M$. 
\medskip 
 
  LRA via random sampling with oversampling  fits this outline.


Apart from the stage of  multiplication by  sketch matrices, performed at a dominated cost,  the algorithm HMT uses 
 $(k+p)m$ scalars and $O((k+p)^2m)$
 flops, while                                                               the TYUC and GN algorithms  use about
 $mk+n\ell$ and  $mk+n\ell+1.5kl$  scalars, respectively, and $O((m+k)k\ell)$ 
 flops, for $k\le \ell$.  
Unlike the number of scalars, the bound on the number of flops decreases 
greatly where $k+p\ll n$ and $k\le \ell \ll n$ for HMT and for TYUC and GN algorithms,  respectively. Here and hereafter $a\ll b$ means that $b$ exceeds $a$  greatly in context, e.g., $a = o(b)$ as $b\rightarrow \infty$. In particular the algorithms HMT and  TYUC/GN run at {\em sublinear cost}, that is,  use much fewer than $mn$ 
scalars and flops,
where we seek {\em Very Low Rank Approximation (VLRA)}, namely, where $k+p\ll \sqrt n$ or $kl\ll \sqrt n$, respectively.
\medskip
  
   Alg. \ref{alglratpsvd} 
of Appendix \ref{slrasvd} performs Stage (2) by using about as many
scalars and
$O(m\ell^2)$ flops. We can decrease this bound to $O(mk\ell)$ (which is sublinear for
$kl\ll \sqrt n$)
by applying TUYC or GN algorithms to compute near-optimal rank-$r$ approximation of $M(\rho)$ rather than its $r$-top SVD at stage (2)
of Alg. \ref{algesc1}.
 
\subsection{Accuracy of the output LRA}\label{saccrc}

 By virtue of Thm. \ref{thcrr2sp1},
the spectral norm 
$||\cdot||:=||\cdot||_2$ of the 
 output error matrix of Alg. \ref{algesc1}
 satisfies 
\begin{equation}\label{eqerrnrm}
||M-M(\rho)_r||\le ||M-M_r||+2||M-M(\rho)||,
\end{equation}
that is, Alg. \ref{algesc1}
is optimal
up to a factor of
$1+ 2\frac{||M-M(\rho)||}{||M-M_r||}$. 
For $(r,\rho)$-matrices, this factor can be large  only if the  computation of a  reasonably close rank-$\rho$ approximation of $M$ at stage 1 fails. It does fail for worst case matrices $M$
 (see Appendix \ref{shrdin}); moreover,  empirically the output matrices of the algorithms HMT, TYUC, and GN with sparse  sketch matrices tend to be unsatisfactory 
  as  LRA \cite{L09},
  \cite[Sec. 3.9]{TYUC17}, but in our numerical tests in Sec. \ref{ststs}
  accuracy of the output matrices of stage (1) of Alg. \ref{algesc1}
  tended to be  raised  at 
   stage  (2) to near-optimal level for a large sub-class of $(r,\rho)$-matrices $M$.


\subsection{Iterative LRA}\label{sitrlr}

In Sec. \ref{sitref} we  apply our LRA algorithm (or we can apply any efficient LRA algorithm instead) to the output error matrix
of a crude LRA computed at stage (1),  to decrease its norm.
Recursively we  extend
to LRA the popular techniques of iterative refinement  and approximate $M$ with a sum of computed LRAs.
 We  reduce the algorithm essentially to
recursive applications of Alg. \ref{algesc1}, but the  rank of computed sums and hence the  cost of their computation  grow fast as recursive steps advance. We counter this deficiency by recursively compressing the sums into rank-$r$ matrices, although this  limits
 output accuracy of computed LRAs. 
 The resulting iterative LRA algorithm is a natural extension to LRA of an efficient popular technique of matrix computations; we specify, analyze, and  test it in Secs. \ref{sitref} and  \ref{ststs}.

Our upper estimates for the accuracy of both Alg. \ref{algesc1} and iterative refinement are close to one another, although both can be quite pessimistic (see Remarks \ref{resrp} and  \ref{resrpk}).  

 In our tests, the output accuracy of our iterative LRA  applied to $(r,\rho)$-matrices 
improved significantly 
at the second iteration,
 little improved                                                                                                                                                                                                                                                     afterwards, 
and overall error norm  was never less  than with  Alg. \ref{algesc1}  for   
the upper rank $\rho=4r$. 

\subsection{Related works}\label{srelwrk}
    

 Our iterative  LRA  extends the algorithms of \cite{LP19,PL19}, which in turn extend to LRA the
classical  technique
of iterative refinement, particularly popular 
 for the solution of a linear system of equations (see 
\cite[Secs. 3.3.4 and 4.2.5]{S98}, 
\cite[Ch. 12 and Sec. 20.5]{H02}, \cite[Secs. 3.5.3 and 5.3.8]{GL13}) but also widely applied 
 in various other linear and nonlinear matrix computations (see \cite[page 223 and the references on page 225]{S98}). 
 Its extension to LRA is not  
 straightforward  
 because one must apply {\em re-compression} to block recursive increase of the rank of the computed approximations; in a sense we extend re-compression of \cite{P93,PRW02}.
Such compression is commonly 
used for LRA with oversampling. 
For sparse test matrices we  used
abridged SRHT matrices  studied
 in the papers  \cite{PLSZ16,PLSZ17,
PLSZa}, which in turn extends the earlier  study of
 Gaussian elimination without pivoting
 in \cite{PQY15,PZ17}.  

 Unlike our presentation
in Secs.
\ref{sprgr} and \ref{sitrlr}
as well as in Secs.
\ref{sspdsc} and \ref{sitref}, we actually first  devised
 iterative refinement of LRA that incorporated our oversampling technique and then noticed that that technique alone can be  competitive or  superior.
 
 
\subsection{Organization  
of our paper}
 We present our accelerated oversampling LRA algorithm in Sec. \ref{sspdsc},   
discuss fast computation of a crude LRA in Sec. \ref{sspfcrd},  present our iterative LRA algorithm in Sec. \ref{sitref},
and cover our numerical experiments in Sec. \ref{ststs}. We devote our short Sec. \ref{scncl} to conclusions.
In  Appendix \ref{shrdin}
 we describe some small families of  matrices for which any
 LRA algorithm  fails unless it access all entries of an input matrix. In  Appendix \ref{slrasvd} we fast compute $r$-top SVD of an LRA of a low-rank matrix. 
 In Appendix \ref{spreprmlt}
   we generate abridged
   SRHT matrices, whose incorporation as  sketch matrices makes random sampling LRA  run faster. 


\section{LRAs of $(r,\rho)$-matrices: generic algorithm and its analysis}\label{sspdsc}

We begin with some  background materials.
  
  Hereafter 
 $||X||:=\sigma_1(X)$ denotes the spectral norm  of a matrix $X$, to which we restrict our study, except for Remark \ref{rerho1}.
In particular we only use  Eckart-Young's (rather than Mirski's) theorem:                                                                                                                                                                                                                                                                                                                                                                                                                                                                                                                                                                           

\begin{theorem}\label{theym}
$\sigma_{r+1}(M)=||M_{r}-M||=
\min_{X:~\rank(X)\le r} ||M-X||$.
\end{theorem}


 \begin{theorem}\label{thsngr} {\rm \cite[Cor. 8.6.2]{GL13}.}
 For a pair of $m\times n$ matrices $M$ and $M+E$ it holds that
 $$|\sigma_j(M+E)-\sigma_j(M)|\le||E||~{\rm for}~j=1,\dots,\min\{m,n\}. $$
  \end{theorem} 

In the following algorithm we elaborate upon the outline of Sec. \ref{sprgr} assuming that we are given  black box subroutines
 that output a crude LRA  and an estimate for a spectral norm
 of a matrix.

\begin{algorithm}
\label{algesc1} LRA  of $(r,\rho)$-matrices.
 

\begin{description}


\item[{\sc Input:}] 
An $m\times n$ matrix $M$, for $n\le m$,  given implicitly by an oracle (black box subroutine) for its multiplication by a vector (cf. Sec. 1.1.1), a target rank $r\ll n$, and a tolerance 
$\mu$ to the output error norm.


\item[{\sc Output:}]
FAILURE or  SVD of a rank-$r$
matrix $X\approx M$ such that $||M-X||\le\mu$.


\item[{\sc Initialization:}]
Fix an 
{\rm upper rank}  
 $\rho$ satisfying   
\begin{equation}\label{eqtau}
r\le \rho\ll \sqrt n.
\end{equation}


\item[{\sc Computations:}]

\begin{enumerate}
\item
Compute  an approximation of $M$ given by SVD of a matrix  $M(\rho)$ of   rank $\rho$.
\item  Compute its $r$-top SVD (cf. 
\cite[Appendix B]{PL19}), which is SVD of $X:=M(\rho)_r$.
\item  
Estimate the  norm $||E||$. If
\begin{equation}\label{eqmu}
||E||\le \mu~{\rm for}~
E:=M-X=M-M(\rho)_r,
\end{equation}
then output 
that SVD 
and $||E||$.
 Otherwise
output FAILURE.
\end{enumerate}
\end{description}
\end{algorithm}

Next we  prove bound (\ref{eqerrnrm}).    
\begin{theorem}\label{thcrr2sp1}
The error norm  $||E||=||M-X||$ (cf.~Eqn. (\ref{eqmu})) satisfies
\begin{equation}\label{eqrlerrnrm}
||E||\le \sigma_{r+1}(M)+2\nu~{\rm where}~ \nu=\nu(M,\rho):=||M-M(\rho)||,
\end{equation}
which turns into 
(\ref{eqerrnrm})
for $\sigma_{r+1}(M)=||M-M_r||$.

\end{theorem}


 \begin{proof}   
Represent $E$ of (\ref{eqmu}) as $M-M(\rho)+M(\rho)-M(\rho)_r$,
apply triangle inequality, and obtain 
$$||E||\le ||M-M(\rho)||+||M(\rho)-M(\rho)_r||=\nu+\sigma_{r+1}(M(\rho)).$$
 
Thm. \ref{thsngr} implies that 
$\sigma_{r+1}(M(\rho))\le \nu+\sigma_{r+1}(M)$.
Combine the latter two bounds. 
\end{proof}


\begin{remark}\label{resrp}
Upper bound (\ref{eqrlerrnrm})
is within $2 \nu$ from the lower
bound  $\sigma_{r+1}(M)$ but can be quite pessimistic,  e.g.,  $||E||=\sigma_{r+1}(M)$  where $M(\rho)=
M_{\rho}$.
\end{remark}

\begin{remark}\label{rerho1}
{\rm [The  choice of an upper rank $\rho$.]}  We should choose  a  larger 
value $\rho$ to decrease the  norm bound $\nu$ but should choose it smaller
to perform  stage 2 faster. For a try-and-error selection among various candidate values of
$\rho$, we can  compare either the ratios   $\frac{||M-M(\rho)||_1}{||M-M(r)||_1}$
or, at a higher computational cost, the values
 $\nu(M,\rho)$.
\end{remark}

\begin{remark}\label{rerho} {\em [A simplified version.]}   In Sec. \ref{stest1} we
 test a simplified version of Alg. \ref{algesc1} where we
  impose no upper bound  $\mu$  and  output  the error norms $||E||=||E(\rho)||$ for $\rho=ir$ and $i=2,3,4,5$. In practice, the choice of  integers $\rho$ and then $r$ should be adapted
  towards minimizing 
the error norms of the computed rank-$\rho$ and rank-$r$
approximations of $M$.
\end{remark}  

\section{LRA
by means of random sampling with sparse  sketch matrices}
\label{sspfcrd}

                                                                                                                                                                                                                                                                                                                                                                                                                                                                                                                                                                                                                                                                                                                                                                                                                                                                                                                                                                                                                                                                                                                                                                                                                                                                                                                                                                                                                                                                                                                                                                                                                                                                                                                                                                                                                                    The random sampling LRA algorithms of \cite{HMT11,TYUC17} are close to optimal with a high probability under various  customary choices of dense sketch matrices such as  Gaussian, Rademacher, SRHT,  and 
SRFT. In the case of VLRA (where $(r+p)^2\ll n$ in the HMT algorithm or $kl\ll n$ in the TUYC and GN algorithms)
one can accelerate the computations 
 by using sparse  sketch matrices, such as Ultra-Sparse Rademacher  sketch matrices \cite[Sec. 3.9]{TYUC17}, \cite{CW09}, which
  depend on a positive integer  parameter $s$, such that sketch matrices are sparse and support  faster LRA for smaller values $s$.
  In the present paper we tested just TYUC algorithm with  sparse {\em abridged}   SRHT  matrices, specified in Appendix \ref{spreprmlt} and also supporting faster VLRA. For other candidate sparse  sketch matrices 
see \cite[Remark 4.6]{HMT11}, \cite[Sec. 5.1.8]{GL13}.

According to \cite{L09},
 such acceleration tends to
 make the output LRAs somewhat less reliable, and  in our tests
 LRAs computed at stage 1  of Alg. \ref{algesc1} 
tended to be  cruder  with
abridged SRHT  sketch matrices than with Gaussian random ones. Stage 2, however, tended to fix
this discrepancy where the ratio $\frac{\sigma_{r+1}}{\sigma_{\rho+1}}$ was large enough.
  
  
\section{Iterative  refinement of LRA
of an $(r,\rho)$-matrix  with compression}\label{sitref}

Next we elaborate 
upon iterative refinement of LRA of an $(r,\rho)$-matrix, sketched in Sec. \ref{sitrlr}.
As a basis we apply 
Alg. \ref{algesc1}, but other efficient LRA algorithms can be applied instead of it.  We  allow at most $h$ refinement steps for a fixed tolerance $h$. 
 As usual for iterative refinement we perform some additions/subtractions at updating steps  with a higher precision, and we increase a lower bound on the upper rank $\rho$ from $r+1$ to $2r$. 

\begin{algorithm}
\label{algesck} Iterative  refinement of LRA of an $(r,\rho)$-matrix.
 

\begin{description}


\item[{\sc Input and OUTPUT}]
as in Alg. \ref{algesc1}, except that the input of  Alg. \ref{algesck} additionally includes a positive integer $h$, which bounds the number of recursive steps allowed, and an Alg. $\mathbb S$ (any black box subroutine) that with a higher precision adds or subtracts a pair of black box $m\times n$ matrices.


\item[{\sc Initialization:}]

Write $i:=0$, $X_0:=0$,  and $E_0:=M$ and fix an  upper rank $\rho$ 
satisfying
\begin{equation}\label{eqtauk}
2r\le\rho\ll n.
\end{equation}


\item[{\sc Computations:}]
\begin{enumerate}
\item
Compute an approximation of $E_i$  by  SVD of a matrix  $E_i(\rho)$ of  rank $\rho$. [This is step 1 of Alg. \ref{algesc1} applied to $E_i$ rather than $M$.]

\item
Compute the matrix 
$M_i(\rho):=X_i+E_i(\rho)$
by applying Alg. $\mathbb S$, performed with a higher precision. [This updates approximation $X_i$ to $M$.]
\item
Compute  SVD of  $X_{i+1}:=
(M_i(\rho))_r$. [This enables us to compress $M_i(\rho)$.] 

\item  Output 
this SVD 
as an LRA of $M$ and also output the norm 
$\mu_i:=||M-(M_i(\rho))_r||$ if it satisfies the tolerance bound, that is, if the  $\mu_i:=\mu_i(M,\rho)\le
\mu$. 
 
Otherwise output FAILURE if $i=h$.
 
\item
Otherwise, update the error matrix
$E_{i+1}:=M-X_{i+1}$, increase
 $i$ by 1, and go to stage 1 (to resume the iteration).
\end{enumerate}

\end{description}
\end{algorithm}

\begin{remark}\label{re01} {\em [The computational cost bound, the choice of an upper rank, and variations  of Alg. \ref{algesck}.]}    As in Alg. \ref{algesc1}  maximize the upper rank $\rho$ to decrease the output error norm, being only limited by the assumption 
    that the  overall computational cost is dominated at stage 1 (cf. (\ref{eqtauk})).  We can allow  an upper rank $\rho$ to vary as $i$ varies as long that it stays 
  much lower than $n$ for all $i$.   In Sec.  \ref{stest1} we test a   
  variant of
Alg.  \ref{algesck} where we  fix $h=3$, and                                                                                                                                                                                                                      
  let $\rho:=r$ at the first iteration (for $i=0$) and 
   $\rho:=2r$ at the next  iterations.
\end{remark}
 
In the following
  extension of Thm. \ref{thcrr2sp1} we  estimate the {\em output error norm} of
Alg. \ref{algesck}.
  
  \begin{theorem}\label{thcrr2spi}
For the error matrices  $E_i:=M-M_{i}(\rho)$ of
Alg. \ref{algesck},
write 
$\nu_i:=||E_i||,~i=0,1,\dots,h-1.$ 
 Then it holds that
\begin{equation}\label{eqnui}
\mu_{i+1}:=||E_{i+1}||\le \sigma_{r+1}(M)+2\nu_i.
\end{equation}
\end{theorem}

\begin{proof}
Recall from  Alg. \ref{algesck} that
$$E_{i+1}:=M-X_{i+1}~{\rm for}~X_{i+1}:=(X_i+
E_i(\rho))_r.$$ 
Hence
$$||E_{i+1}||\le ||M-M_i(\rho)||+||M_i(\rho)-X_{i+1}||=\nu_{i}+\sigma_{r+1}(M_i(\rho)).$$

Thm. \ref{thsngr} implies  that 
$$\sigma_{r+1}(M_i(\rho))\le \sigma_{r+1}(M)+\nu_i.$$
Combine the latter two bounds. 
\end{proof}

\begin{remark}\label{resrpk}
The upper bound  (\ref{eqnui}) on the error norm is within $2 \nu_i$ from the lower bound  $\sigma_{r+1}(M)$,  attained where $M_i(\rho)=M_{\rho}$.
\end{remark}

Next we specify how compression limits the accuracy of the output LRAs.

\begin{lemma}\label{le2}
For any $m\times n$ 
matrix $X$ and $m\ge n$ 
 write $\sigma_j(X):=\infty$ for $j<1$ and $\sigma_j(X):=0$ for $j>n$. Then
$$\sigma_{j+r+1}(M)\le \sigma_{j+1}(E_i)\le 
\sigma_{j-r+1}(M)~{\rm for}~
E_i~{\rm of~Alg.~ \ref{algesck}~and~all}
~i~{\rm and}~j.$$
\end{lemma}

\begin{proof}
Recall that 
$$E_i=M-X_i~{\rm where}~\rank(X_i)\le r~{\rm for~all}~i.$$ Furthermore,
$$\sigma_{j+1}(E_i)=\min_{Y:~\rank(Y)\le j}||E_i-Y||$$
 by virtue of Thm. \ref{theym}.  

Substitute $E_i-Y=M-(Y+X_i)$ where $\rank(Y+X_i)\le j+r$ if $\rank(Y)\le j$ and obtain 
$$\sigma_{j+1}(E_i)\ge\min_{Z:~\rank(Z)\le j+r}||M-Z||=\sigma_{j+r+1}(M).$$ 
Similarly deduce that
$\sigma_{j+1}(E_i)\le 
\sigma_{j-r+1}(M)$.
\end{proof}

\begin{corollary}\label{cocrr2}
The singular values $\sigma_{\rho+1}(E_i)$  satisfy
\begin{equation}\label{eqenrmk}
\sigma_{\rho+r+1}(M)\le \sigma_{\rho+1}(E_i)\le
\sigma_{\rho-r+1}(M),~i=0,1,\dots,h-1.
\end{equation} 
\end{corollary}

 \begin{remark}\label{realh01} {\rm [Error norm bounds.]}  In good accordance with the corollary, in our tests with $(r,\rho)$-matrices
 the output accuracy    of Alg. \ref{algesc1}
 was  higher than  that of  Alg. \ref{algesck} provided that the upper rank $\rho$ in Alg. \ref{algesc1}  was
twice as large as in Alg. \ref{algesck}.
 The test results were mixed for input matrices 
 with flat spectra.
\end{remark}


                                                                                                                                                                                                                                                                                                                                                                                                                                                          
\section{Numerical experiments} \label{ststs}
  
In this section  we cover the tests of our LRA algorithms on  synthetic and real-world  matrices. 
We implemented the algorithms in Python, running them on a 64bit MacOS machine with 16GB Memory and called scipy.linalg for numerical linear algebra routines such as QR factorization with pivoting, Moore-Penrose matrix inversion, and SVD.


\subsection{Input Matrices} \label{sec_input}

\noindent {\bf Real-world Matrices:} 
We computed rank-$r$
approximation for the matrices 
 {\em Shaw} and {\em Gravity} that represent discretization of Integral Equations involved in the built-in problems of the Regularization Tools\footnote{See 
 http://www.math.sjsu.edu/singular/matrices and 
  http://www.imm.dtu.dk/$\sim$pcha/Regutools
  
 For more details see Ch. 4 of 
  http://www.imm.dtu.dk/$\sim$pcha/Regutools/RTv4manual.pdf }. The matrix Shaw is from a one-dimensional image restoration model problem;  the matrix
  Gravity is from a one-dimensional gravity surveying model problem. 
  Both  are 
  $1000\times 1000$ dense real matrices  
having low numerical rank.     
 
Our third 
input matrix, from the discretization of a single layer potential ({\em SLP}) operator, has size  $1024\times 1024$ (see   \cite[Sec. 7.1]{HMT11} for more details).  
  
\begin{figure}
 \includegraphics[width = \textwidth]{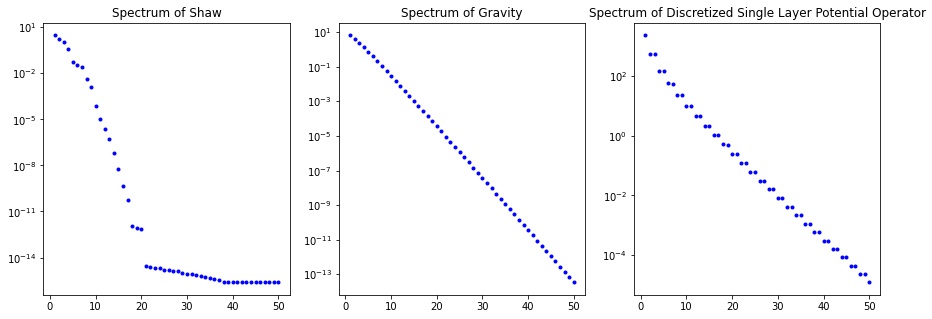}
 \caption{Spectra of 
 singular values of some real-world matrices}
 \label{real_world_spectrum} 
 \end{figure}

Fig. \ref{real_world_spectrum} shows the distribution of the top 50 singular values of these three matrices.

We set the target ranks $r$ equal to 20
for Shaw, 45
for Gravity,  and 11 for SLP matrices.
  
We also computed rank-$r$  approximations for the three matrices called {\em StreamVel},  {\em MinTemp}, and  {\em MaxCut} and used in numerical experiments in \cite[Sec.~7.3.2]{TYUC19}. 

The  $10738 \times 5001$  matrix {\em StreamVel} is involved in numerical simulation (DNS) on a coarse mesh of the 2D Navier–Stokes equations for a low-Reynolds number flow around a cylinder.

The $19264 \times 7305$ matrix {\em MinTemp}  tabulates meteorological variables at weather stations across the northeastern United States on days during the years 1981–2016. 

The $2000\times 2000$ matrix {\em MaxCut}  provides a high-accuracy solution to the MaxCut SDP for a sparse graph \cite{GW95}. 


We let the target rank $r$  be 20 for StreamVel, 
 10 for MinTemp, and 14 for MaxCut matrices.

In our tests, we padded all matrices with 0s to  increase their dimensions to the powers of 2.
\medskip

\noindent {\bf Synthetic Input Matrices.} 

We generated random synthetic $1024\times 1024$ input
matrices  of five classes, two of them  -- with fast and  slowly decaying spectra --   
 as the products $U\Sigma V^T$, 
with the matrices  $U$ and $V$ of the left and right singular vectors of a  Gaussian random matrix, respectively. 
 
By letting
$\Sigma = \textrm{diag}(v)$, $v_i = 1$ for $i = 1, 2, 3,\dots, 20$, 
$v_i = \frac{1}{2^{i-20}}$ for $i = 21, \dots, 100$, and $v_i = 0$ for $i > 100$,
we arrived at the matrices with {\em fast decaying spectra}. 

By letting $\Sigma = \textrm{diag}(u)$, $u_i = 1$ for $i = 1, 2,\dots, 20$, 
and $u_i = \frac{1}{(1 + i - 20)^2}$ for $i > 20$, we 
arrived at the matrices with  {\em  slowly decaying spectra}.

For both input classes we fixed  target rank $r = 20$ and set it equal to 10 for the following three classes  of synthetic $1024\times 1024$ matrices, generated  according to the recipes in  \cite[Sec. 7.3.1]{TYUC19}, which involved
{\em effective rank parameter} $R=20$.

 (i) {\em Low-rank + noise} matrices 
${\bf A} = \diag(\underbrace{1, \dots, 1}_{R}, 0,\dots, 0) + (\xi/1024){\bf C},$
for ${\bf C} = {\bf G}{\bf G}^T$ and a $1024\times 1024$ standard normal  matrix  ${\bf G}$,   depend on a parameter $\xi$  
 and become
{\em Low Rank Low Noise} matrix for $\xi=10^{-4}$, 
 the {\em Low Rank Med Noise} matrix for $\xi=10^{-2}$, and  
 the {\em Low Rank High Noise} matrix for $\xi=10^{-1}$.

(ii)  {\em Polynomial Decay} matrices
${\bf A} = \diag(\underbrace{1, \dots, 1}_{R}, 2^{-p}, 3^{-p}, \dots, (n-R+1)^{-p})$
are  {\em Poly Decay \\Slow} for $p=0.5$,  {\em Poly Decay Med}  for  $p=1$,  
and {\em Poly Decay Fast}
for  $p=2$.

(iii) {Exponential Decay} matrices 
${\bf A} = \diag(\underbrace{1, \dots, 1}_{R}, 10^{-q}, 10^{-2q}, \dots, 10^{-(n-R)q})$
are {\em Exp Decay \\Slow} for $q=0.01$, {\em Exp Decay Med} for $q=0.1$,    and {\em Exp Decay Fast} for $q=0.5$.


\subsection{Two options at stages 1 of  Algs. \ref{algesc1}  and  \ref{algesck}}

In our tests of   Algs.  \ref{algesc1} and  \ref{algesck} we used the following two options  at their stage 1  
for computations  of a crude rank-$\rho$ approximation of a matrix $M$.
  
  (i) Fast  randomized algorithm of \cite{TYUC17} with Gaussian random  sketch matrices
 $F$ of size $(2\rho)\times m$ and $H$ of size $n\times \rho$.
  
  (ii) Its accelerated  randomized variant where                                                                                                                                                                              
 abridged (length 3) SRHT  sketch matrices $H\in\mathbb C^{1024\times \rho}$  replaced Gaussian ones.


\subsection{Test results for Algs. \ref{algesc1} and \ref{algesck}} \label{stest1}

For {\bf Alg. \ref{algesc1}} we output the relative error norms $\text{err} := \frac{||M-(M(\rho))_r||}{||M-M_r||}$
 for  input matrices $M$, a target rank $r$, the $r$-truncation $M_r$  and $M(\rho)_r$,  and  the upper rank $\rho$, set to $2r,3r,4r$, $5r$ (cf. Remark \ref{rerho}).

For {\bf Alg. \ref{algesck}} we output the relative error norms
$\text{err=err}_i: = \frac{||M - X_i ||}{||M - M_r||}$ 
 for the  input matrices $M$ of Alg. \ref{algesck}, the approximations  $X_i$ output in its $i$-th iteration, target ranks $r$, the $r$-truncation $M_r$, and $i=1,2,3$ (letting $h=3$).
 As we pointed out in Remark \ref{re01}, we fixed $\rho=r$  at the initial step, that is,  for $i=1$, and $\rho=2r$  at the other steps, for $i>1$.
 
 
Our tables display the mean  of the relative error norms observed in these tests. 

For  Alg.  \ref{algesc1} 
 we also display the standard deviation.
 
  For both algorithms we repeated the tests 100 times per test set  
  with independent random choices of  sketch matrices $F$ and $H$,
   that is, per the triple  $\{r, \rho, {\rm TYPE}\}$
   for Alg.  \ref{algesc1}
 or the pair $\{r, {\rm TYPE}\}$ for Alg.  \ref{algesck}  where TYPE stands for an input matrix type.   We, however, ran just 10 tests for the StreamVel and MinTemp matrices, which have significantly larger sizes. 

 In Tables \ref{two_stage_test_result},  \ref{it_ref_test_result},
 \ref{two_stage_test_result_TYUC_rw},   and  \ref{it_ref_test_result_TYUC_rw} we display  
  the output error norms in our tests with both Gaussian random and abridged SRHT  sketch matrices, which  were more or less in the same range for both matrix types.   
All our tests with
the synthetic input matrices from \cite[Sec. 7.3.1]{TYUC19} wit
abridged SRHT  sketch matrices have failed, that is, have output 
large error norms. In Tables \ref{two_stage_test_result_TYUC_syn}
and \ref{it_ref_test_result_TYUC_syn} we 
display the results of these tests but only  with Gaussian random  sketch matrices, where the output LRAs were 
 quite  accurate, except for LRAs of   
 high noise and slow decay matrices. 
  
Alg. \ref{algesc1} failed for the matrix Shaw and target rank $r=20$
with both Gaussian and abridged SRHT  sketch matrices. This could have been predicted because the spectrum of the singular values $\sigma_j(M)$ of these matrices $M$ is flat
for $j>20$ (decreasing roughly from $10^{-13}$
to between $10^{-15}$ and $10^{-16}$); in Table \ref{two_stage_test_result} we display  the test results for the matrix Shaw and target rank $r=19$, which are about as favorable as for Gravity and Slow Decay. 

In both cases of  success and failure
of both algorithms,  {our test results were in good accordance with our formal study}.
 The relative   error norms output by both Algs. \ref{algesc1} and 
\ref{algesck} applied to $(r,\rho)$-matrices have consistently decreased
and frequently to the optimal value 1 (up to  the measurement errors) as  we increased the upper rank $\rho$ for  Alg. \ref{algesc1} and in the transition from the first to the second   iteration of Alg. \ref{algesck}.
In all tests except for High Noise, Slow Decay, and 
 MinTemp matrices already err$_2$ decreased below 1.1.
 The third iteration decreased the error norm much less or even slightly increased it.

For $(r,\rho)$-matrices the error norm bounds of  Alg. \ref{algesc1} for $\rho=4r$ were consistently smaller than those of  Alg. \ref{algesck} in all three iterations.
This was not the case just for the matrices with flat spectra such as  the matrices Shaw for $r=20$, High Noise, and Exp Decay Slow.  

\begin{table}[h!]
\footnotesize
\begin{center}
\begin{tabular}{c| c| l| l| l| l}
\multirow{2}{*}{Input Matrix} 
&  \multirow{2}{*}{Sketch Matrix} 
& \multicolumn{4}{c}{$\rho$}\\  
& 
& \multicolumn{1}{c|}{$2r$} 
& \multicolumn{1}{c|}{$3r$} 
& \multicolumn{1}{c|}{$4r$} 
& \multicolumn{1}{c}{$5r$}\\
\hline \hline
\multirow{2}{*}{Gravity} 
& Abr. SRHT
& 1.000 $\pm$ 1.036E-05
& 1.000 $\pm$ 7.815E-06
& 1.000 $\pm$ 7.888E-06
& 1.000 $\pm$ 6.683E-06\\

& Gaussian

& 1.000 $\pm$ 7.495E-06
& 1.000 $\pm$ 5.211E-06
& 1.000 $\pm$ 7.214E-06
& 1.000 $\pm$ 5.432E-06\\
\hline

\multirow{2}{*}{SLP} 
& Abr. SRHT
& 1.970 $\pm$ 4.083
& 1.000 $\pm$ 0.002
& 1.000 $\pm$ 1.018E-07
& 1.000 $\pm$ 1.362E-10 \\

& Gaussian
& 1.001 $\pm$ 0.0008160
& 1.000 $\pm$ 1.528E-07
& 1.000 $\pm$ 1.075E-10
& 1.000 $\pm$ 6.350E-14\\
\hline

\multirow{2}{*}{Shaw,  $r=19$}
& Abr. SRHT
& 1.000 $\pm$ 1.778E-06
& 1.000 $\pm$ 1.359E-06
& 1.000 $\pm$ 1.665E-06
& 1.000 $\pm$ 1.441E-06\\

& Gaussian
& 1.000 $\pm$ 2.256E-06
& 1.000 $\pm$ 1.685E-06
& 1.000 $\pm$ 1.260E-06
& 1.000 $\pm$ 1.427E-06\\

\hline
\multirow{2}{*}{Fast Decay} 
& Abr. SRHT
& 1.000 $\pm$ 3.086E-12
& 1.000 $\pm$ 3.590E-16
& 1.000 $\pm$ 2.470E-16
& 1.000 $\pm$ 2.617E-16\\
& Gaussian
& 1.000 $\pm$ 6.651E-12
& 1.000 $\pm$ 5.117E-16
& 1.000 $\pm$ 4.552E-16
& 1.000 $\pm$ 4.590E-16\\
\hline
\multirow{2}{*}{Slow Decay} 
& Abr. SRHT
& 1.000 $\pm$ 3.604E-05
& 1.000 $\pm$ 1.244E-06
& 1.000 $\pm$ 1.543E-07
& 1.000 $\pm$ 3.754E-08\\
& Gaussian
& 1.000 $\pm$ 4.130E-05
& 1.000 $\pm$ 1.380E-06
& 1.000 $\pm$ 1.783E-07
& 1.000 $\pm$ 3.708E-08\\

\end{tabular}
\end{center}
\caption{Test results 
for Alg. \ref{algesc1} on Gravity,  
discretized single layer potential operator (SLP), Shaw, and
synthetic matrices Fast Decay and Slow Decay, with fast and 
slowly decaying spectra.}
\label{two_stage_test_result}
\end{table}

\begin{table}[h]
\footnotesize
\begin{center}
\begin{tabular}{c| c | c | c | c}
Input Matrix &  Sketch Matrix & 1st Itr. & 2nd Itr. & 3rd Itr. \\
\hline\hline
\multirow{2}{6em}{Fast Decay} 
& Abr. SRHT & 3.155000 & 1.0000 & 1.0000 \\
& Gaussian & 3.1202 & 1.0000 & 1.0000 \\
\hline
\multirow{2}{6em}{Slow Decay} 
& Abr. SRHT & 5.0468 & 1.0003 & 1.0001 \\
& Gaussian & 5.0755 & 1.0002 & 1.0001 \\
\hline
\multirow{2}{4em}{Shaw for $r=20$} 
& Abr. SRHT & 2.8820E+01 & 1.0983 & 1.1225 \\
& Gaussian & 1.8235E+01 & 1.1517 & 1.1189 \\
\hline
\multirow{2}{4em}{Gravity} 
& Abr. SRHT & 1.5762E+01 & 1.0000 & 1.0000 \\
& Gaussian & 1.2917E+01 & 1.0000 & 1.0000 \\
\hline
\multirow{2}{4em}{SLP} 
& Abr. SRHT & 1.0931E+02 & 1.0014 & 1.0000 \\
& Gaussian & 5.2205 & 1.0000 & 1.0000 \\

\end{tabular}
\caption{Test results  
for Alg.  \ref{algesck}
 on Fast Decay, Slow Decay, Shaw, Gravity, and SLP matrices.}
\label{it_ref_test_result} 
\end{center}
\end{table}

\begin{table}[h!]
\footnotesize
\begin{center}
\begin{tabular}{c| c| l| l| l}

 \multirow{2}{*}{Input Matrix}
& \multicolumn{4}{c}{$\rho$}\\  

& \multicolumn{1}{c|}{$2r$} 
& \multicolumn{1}{c|}{$3r$} 
& \multicolumn{1}{c|}{$4r$} 
& \multicolumn{1}{c}{$5r$}\\

\hline \hline

Low Rank Low Noise
& 1.0416 $\pm$ 8.9977E-02
& 1.0000 $\pm$ 2.1197E-06
& 1.0000 $\pm$ 2.4039E-06
& 1.0000 $\pm$ 2.2834E-06

\\
\hline
Low Rank Med Noise
& 1.4335 $\pm$ 1.7048E-01
& 1.0382 $\pm$ 3.1809E-02
& 1.0057 $\pm$ 1.4442E-03
& 1.0026 $\pm$ 5.8650E-04
\\
\hline
Low Rank High Noise
& 5.6972 $\pm$ 8.6182E-01
& 4.8401 $\pm$ 4.3819E-01
& 4.0328 $\pm$ 2.3493E-01
& 3.7893 $\pm$ 2.1626E-01
\\
\hline
 Poly Decay Slow
& 2.0588 $\pm$ 1.8783E-01
& 1.6525 $\pm$ 2.0027E-01
& 1.3617 $\pm$ 8.6188E-02
& 1.2062 $\pm$ 8.6352E-02
\\
\hline
 Poly Decay Med
& 1.5384 $\pm$ 2.1907E-01
& 1.0315 $\pm$ 2.6585E-02
& 1.0028 $\pm$ 1.0190E-03
& 1.0009 $\pm$ 4.5652E-04
\\
\hline
 Poly Decay Fast
& 1.3133 $\pm$ 1.5431E-01
& 1.0001 $\pm$ 1.3887E-04
& 1.0000 $\pm$ 3.7489E-06
& 1.0000 $\pm$ 2.4495E-07
\\
\hline
 Exp Decay Slow
& 2.8587 $\pm$ 3.3389E-01
& 2.2772 $\pm$ 2.0481E-01
& 1.8244 $\pm$ 1.0970E-01
& 1.5721 $\pm$ 1.0528E-01
\\
\hline
 Exp Decay Med
& 1.5576 $\pm$ 1.2324E-01
& 1.0414 $\pm$ 4.9018E-02
& 1.0001 $\pm$ 9.0188E-05
& 1.0000 $\pm$ 3.7953E-07
\\
\hline
 Exp Decay Fast
& 1.3121 $\pm$ 1.4989E-01
& 1.0000 $\pm$ 6.4663E-11
& 1.0000 $\pm$ 3.6020E-16
& 1.0000 $\pm$ 3.0986E-16
\\

\end{tabular}
\end{center}
\caption{Test results 
for Alg. \ref{algesc1} on  the synthetic input matrices from \cite[Sec. 7.3.1]{TYUC19} using Gaussian  sketch matrices.}
 \label{two_stage_test_result_TYUC_syn}
\end{table}

\begin{table}[h]
\footnotesize 
\begin{center}
\begin{tabular}{c| c | c | c}
Input Matrix & 1st Itr. & 2nd Itr. & 3rd Itr. \\

\hline\hline

Low Rank Low Noise
& 1.3940  & 1.0000 & 1.0000 \\
\hline
Low Rank Med Noise
& 1.4752 & 1.0386 & 1.0375 \\
\hline
Low Rank High Noise
& 1.4507  & 1.5135  & 1.5332 \\
\hline
Poly Decay Slow
& 1.5920  & 1.4154 & 1.4073 \\
\hline
Poly Decay Med
& 1.5569 & 1.0345 & 1.0306 \\
\hline
Poly Decay Fast
& 1.3784  & 1.0001  & 1.0002 \\
\hline
Exp Decay Slow
& 1.5202  & 1.4956 & 1.4750 \\
\hline
Exp Decay Med
& 1.4946 & 1.0115 & 1.0164 \\
\hline
Exp Decay Fast
& 1.4235 & 1.0000 & 1.0000 \\

\end{tabular}
\end{center}
\caption{Test results for Alg. \ref{algesck} on  the synthetic input matrices from \cite[Sec. 7.3.1]{TYUC19} using Gaussian  sketch matrices.}
\label{it_ref_test_result_TYUC_syn}
\end{table}



\begin{table}[h!]
\footnotesize
\begin{center}
\begin{tabular}{c| c| l| l| l| l}
Input 
&  Sketch Matrix
& \multicolumn{4}{c}{$\rho$}\\  
Matrix
&  Type
& \multicolumn{1}{c|}{$2r$} 
& \multicolumn{1}{c|}{$3r$} 
& \multicolumn{1}{c|}{$4r$} 
& \multicolumn{1}{c}{$5r$}\\
\hline \hline
\multirow{2}{*}{StreamVel}
& Abr. SRHT
& 1.5840 $\pm$ 1.2583
& 1.4736 $\pm$ 1.4538
& 1.0022 $\pm$ 2.6478E-03
& 1.0002 $\pm$ 2.7733E-04
\\
& Gaussian
& 1.0932 $\pm$ 5.8887E-02
& 1.0064 $\pm$ 2.9809E-03
& 1.0013 $\pm$ 1.4141E-03
& 1.0003 $\pm$ 3.1762E-04
\\
\hline
\multirow{2}{*}{MinTemp}
& Abr. SRHT
& 2.4575 $\pm$ 2.5748E-01
& 1.6980 $\pm$ 2.3998E-01
& 1.2919 $\pm$ 1.3606E-01
& 1.1181 $\pm$ 9.7415E-02
\\
& Gaussian
& 2.4393 $\pm$ 3.0823E-01
& 1.7216 $\pm$ 2.2510E-01
& 1.3793 $\pm$ 7.6629E-02
& 1.1238 $\pm$ 8.3331E-02
\\
\hline
\multirow{2}{*}{MaxCut}
& Abr. SRHT
& 1.0024 $\pm$ 1.3630E-03
& 1.0005 $\pm$ 1.6875E-04
& 1.0002 $\pm$ 6.8092E-05
& 1.0001 $\pm$ 2.9040E-05
\\
& Gaussian
& 1.0026 $\pm$ 1.3856E-03
& 1.0005 $\pm$ 2.1129E-04
& 1.0002 $\pm$ 6.6260E-05
& 1.0001 $\pm$ 2.4207E-05
\\
\end{tabular}
\end{center}
\caption{Test results 
for Alg. \ref{algesc1} on  the 
real-world input matrices from \cite[Sec. 7.3.2]{TYUC19}}
\label{two_stage_test_result_TYUC_rw}
\end{table}


\begin{table}[h]
\footnotesize 
\begin{center}
\begin{tabular}{c|c| c | c | c}
Input Matrix &  Sketch Matrix Type & 1st Itr. & 2nd Itr. & 3rd Itr. \\
\hline\hline

\multirow{2}{*}{StreamVel} 
& Abr. SRHT 
& 2.9700
& 1.0739 
& 1.0300 
\\
 & Gaussian
& 2.8724
& 1.0397 
& 1.0177 
\\
\hline
\multirow{2}{*}{MinTemp}
& Abr. SRHT
& 2.3340 
& 1.5961 
& 1.5079 
\\
& Gaussian
& 2.3120 
& 1.6410 
& 1.4950 
\\
\hline
\multirow{2}{*}{MaxCut} 
& Abr. SRHT 
& 3.6200E+01 
& 1.0312 
& 1.0182 
\\
& Gaussian 
& 2.6202E+01 
& 1.0302 
& 1.0176\\

\end{tabular}
\end{center}
\caption{Test results  
for Alg.  \ref{algesck}
for $\rho=2r$ on real-world input matrices from \cite[Sec.~7.3.2]{TYUC19}. }
\label{it_ref_test_result_TYUC_rw}
\end{table}




  
\section{Conclusions
}\label{scncl} 

There are various natural extensions
of our study, which can be practically competitive. E.g., one can 
  devise and test
 combinations of Alg.  \ref{algesc1} and of  the TYUC and GN algorithms  with various sparse  sketch matrices such as 
 Ultra-Sparse Rademacher  sketch matrices \cite[Sec. 3.9]{TYUC17}, sparsified random Givens rotations by Rokhlin and Tygert (see \cite[Remark 4.6]{HMT11}), or even further sparsified count sketch matrices in Clarkson and Woodruff \cite{CW09,CW17}.\footnote{\cite{PLSZ16,PLSZ17,PLSZa} generated abridged SRFT matrices similarly to abridged SRHT matrices and in extensive tests  observed similar efficiency of both of these matrix classes.}  Instead of random sampling 
 LRA with  sparse sketch matrices one can apply  Cross-Approximation algorithm for LRA (see \cite{OZ18} and the references therein).

Our comparative study
has showed  no significant benefits of applying  iterative refinement to LRA versus 
  standalone Alg.  \ref{algesc1},
but one can further explore this
comparison. 

\medskip
\medskip
\medskip  
\medskip
\medskip
\medskip


{\bf \Large Appendix} 
\appendix

\section{Small families of matrices that are hard for fast LRA}\label{shrdin}

  Any algorithm that does not access all entries of an input matrix
  fails to compute a close LRA of  the following small families of matrices.
  
\begin{example}\label{exdlt} 
 Let  $\Delta_{i,j}$ denote an $m\times n$ matrix
 of rank 1  filled with 0s except for its $(i,j)$th entry filled with 1. The $mn$ such matrices $\{\Delta_{i,j}\}_{i,j=1}^{m,n}$ form a family of  $\delta$-{\em matrices}.
We also include the $m\times n$ null matrix $O_{m,n}$
filled with 0s  into this family.
Now fix any   algorithm that does not access the $(i,j)$th  
entry of its input matrices  for some pair of $i$ and $j$. Such n algorithm outputs the same approximation 
of the matrices $\Delta_{i,j}$ and $O_{m,n}$,
with an undetected  error at least 1/2.
We arrive at the same conclusion by applying the same argument to the
set of $mn+1$ small-norm perturbations of 
the matrices of the above family and to the                                                                                                                                   
 $mn+1$ sums  
of the latter matrices with  any
   fixed $m\times n$ matrix of low rank. Likewise,
  we can verify that any
 randomized LRA
algorithms that does not access ail entries of an input matrix
fails on those matrix families with 
error probability not close to 0.
\end{example}

\begin{remark}\label{redlt} 
How representative are the above matrix families of hard inputs? The matrices of these families represent data singularities 
and are not relevant to a large class of  input matrices representing regular processes or, say, smooth surfaces. Furthermore, Example   
 \ref{exdlt} does not apply to various important special classes of input matrices such as symmetric or diagonally dominant.  
\end{remark}


\section{Computation of  $r$-top SVD 
 of a matrix product}\label{slrasvd}


For completeness of our exposition we next  recall the  classical algorithm that
 computes $r$-top SVD of a matrix product
 $AB$ (cf. \cite{P93,PRW02})
based on computation of the SVDs of $A$ and $B$.

\begin{algorithm}\label{alglratpsvd}
{\rm [Computation of an $r$-top  SVD of a matrix product.]}
 

\begin{description}


\item[{\sc Input:}]
Four integers $r$, $k$, $m$, and $n$ such that
$0<r\le k\le \min\{m,n\}$ and
two matrices 
$A\in \mathbb C^{m\times k}$ and  
$B\in \mathbb C^{k\times n}$. 

\item[{\sc Output:}]
Three matrices  
$U\in \mathbb C^{m\times r}$ (unitary),  
$\Sigma\in \mathbb C^{r\times r}$ 
(diagonal), and  
$V\in \mathbb C^{n\times r}$  (unitary) such
that $(AB)_{r}=U\Sigma V^*$ is  an $r$-top SVD of $AB$.                                                                                                                                                    


\item[{\sc Computations:}]
\begin{enumerate}
\item
 Compute SVDs $A=U_A\Sigma_AV^*_A$ and $B=U_B\Sigma_BV^*_B$
 where
$U_B\in \mathbb C^{m\times k}$, 
$V_B^*\in \mathbb C^{k\times n}$, and 
 $\Sigma_A,V^*_A,U_B,\Sigma_B\in \mathbb C^{k\times k}$.
   \item
Compute $k\times k$ matrices $W=\Sigma_AV^*_AU_B\Sigma_B$, 
$U_W$, $\Sigma_W$, and $V^*_W$
such that $W=U_W\Sigma_WV_W^*$ is SVD, 
$\Sigma_W=\diag(\Sigma,\Sigma')$, and
 $\Sigma=\diag(\sigma_j)_{j=1}^{r}$ and 
$\Sigma'=\diag(\sigma_j)_{j=r+1}^{k}$
are the matrices of the $r$ top (largest) and the $k-r$ 
trailing (smallest)
singular values of the matrix $W$, respectively.
Output the matrix $\Sigma$.
 \item
Compute and output the matrices $U$ and $V$ made up of the first $r$ 
columns of the matrices $U_AU_W$ and $V_BV_W$, respectively.
\end{enumerate}

\end{description}
\end{algorithm}


The algorithm involves 
$(m+n+k)k$ scalars and 
$O((m+n)k^2)$ flops
(cf. \cite[Figure 8.6.1]{GL13}), running 
particularly fast where $k\ll n\le m$.
Its correctness follows from equations $AB=U_AWV^*_B$, $W=U_W\Sigma_W V^*_W$,
and $\Sigma_W=\diag(\Sigma,\Sigma')$.   
 One can readily extend this algorithm to computing $r$-top SVD of a matrix product
 $ABC$  by using $O(mk^2+n\ell^2)$ flops
 and
$(m+k)k+(n+\ell)\ell$ scalars  where $A\in \mathbb  R^{m\times k}$, $B\in \mathbb  R^{k\times \ell}$, and $C\in \mathbb  R^{\ell\times n}$.

                                                                                                                                                                                                                                                                                                                                                                                                                                                                                                                                                                                                                                                                                                                                                                                                                                                                                                              
\section{Generation of abridged SRHT  matrices}\label{spreprmlt}


In this section  we specify the family of 
abridged SRHT matrices, used in our tests for Algs. \ref{algesc1} and \ref{algesck}. They are defined 
  by  means of abridging 
    the classical 
recursive processes of the generation  of $n\times n$ 
SRHT  matrices,
  obtained from   the  $n\times n$ dense
matrices $H_n$ of  Walsh-Hadamard transform for $n=2^t$
(cf. \cite[Sec. 3.1]{M11}).
The $n\times n$ matrices $H_n$ are obtained in $t=\log_2(n)$ recursive  steps, but we 
 only perform $d\ll t$
 steps, and the resulting 
 abridged matrix $H_{d,d}$ can be multiplied by a vector by using $2dn$   additions and subtractions.
 SRHT matrices are obtained
 from the matrices $H_n$ by means of random sampling and scaling, which we also apply to the $d$-{\em abridged Hadamard  transform matrices}  $H_{d,d}$. They turn into $H_n$ for $d=t$ but are sparse fo $d\ll t$. Namely, we
   write $H_{d,0}:=I_{n/2^d}$ and then
specify the following recursive process: 
 

\begin{equation}\label{eqrfd}
H_{d,0}:=I_{n/2^d},~
H_{d,i+1}:=\begin{pmatrix}
H_{d,i} & H_{d,i} \\
H_{d,i} & -H_{d,i}
  \end{pmatrix}
  ~{\rm for}~i=0,1,\dots,d-1, 
\end{equation}

For any fixed pair of $d$ and $i$, 
each of the matrices  
 $H_{d,i}$ 
is orthogonal 
up to scaling and
 has $2^d$ nonzero entries 
in every row and  column;
  we can compute the product 
 $MH$ for an   $n\times k$ submatrix $H$ of $H_{d,d}$ by using less than $km2^d$ additions and subtractions.


Now define the 
$d$-Abridged Scaled and Permuted 
 Hadamard matrices, $PDH_{d,d}$, where $P$ is a
 random 
 sampling 
 matrix and  $D$ is the matrix of  
random integer diagonal scaling.
Each random permutation or scaling  
 contributes up to $n$ random parameters.  
We can involve more random parameters by applying random permutation and scaling 
also to some or all
 intermediate matrices $H_{d,i}$ for $i=0,1,\dots,d$.

The first $k$ columns of $H_{d,d}$  for  
 $r\le k\le n$ form a
 $d$-Abridged Subsampled Randomized Hadamard Transform (SRHT) matrix $H$, which turns into
a SRHT matrix for $d=t$, where $k=r+p$, $r$ is a target rank and $p$ is the oversampling parameter
 (cf. \cite[Sec. 11]{HMT11}). 

  \medskip
 
{\bf  Acknowledgement:} Jianlin Xia helped place our study into proper historical context.
 




\end{document}